# On an Unnoticed Geometrical Paradox in *Timaeus*' Cosmology

**Luc Brisson**
Centre Jean Pépin
CNRS-ENS/UMR 8230

**Salomon Ofman**
Institut mathématique de Jussieu-
Paris Rive Gauche/HSM
Sorbonne Université

## Abstract

In a previous article, we discussed a paradox in *Timaeus*' cosmology: the fact that there is no void *inside* the universe, even though it is entirely filled with polyhedra—which is a mathematical impossibility (Brisson and Ofman 2025). In the present article, we examine another paradox. While the first paradox is well known and was already highlighted by Aristotle as a fundamental mathematical contradiction undermining Plato's cosmology, this new paradox has gone almost entirely unnoticed by commentators, both ancient and modern. This oversight may surprise scholars, given the extensive body of work on *Timaeus*' universe, much of which emphasizes discrepancies with astronomical observations or points out alleged internal contradictions. Like the first paradox, this one arises from the premise of a universe entirely filled with polyhedra. However, in this case, the contradiction stems from the absence of void *outside* it.

In the first section, we demonstrate that the shape of the universe cannot be a perfect mathematical sphere: that is, its boundary is not smooth but exhibits bumps and hollows. Next, we present conceptual arguments from Plato's text that support the necessity of such 'defects' in the universe's shape as compared to a perfect mathematical sphere. In the third section, we argue that such a universe cannot move at all. Finally, we propose a solution to this mathematical contradiction in *Timaeus*' construction, drawing on the same ideas used to address the earlier apparent contradiction: the unique feature of *Timaeus*' universe as a living being, whose parts are continuously moving, changing, decomposing, and reforming.

While this problem does not depend on the various trends of interpretation of the *Timaeus*, it is related to some important issues concerning Plato's philosophy. These issues include the importance of observations in science—particularly in astronomy—the relationship between intelligible models and their sensible copies, the *mythos/ logos* approach of Plato's cosmology, and the debate over 'metaphorical' *vs* 'literal' interpretation. Of course, all these questions fall outside the scope of this article and will not be addressed here.

## Introduction

In this article, we aim to demonstrate that two fundamental features of *Timaeus*' universe—on the one hand, the absence of void, on the other hand, the complete filling of the universe by polyhedral particles—lead to a mathematical impossibility when considered alongside another feature: the rotation of the universe. In fact, the first two features already contradict another characteristic of the universe: its spherical shape. Both paradoxes originate from two assertions: the four regular polyhedra composing all the bodies of the universe fill it entirely, and the absence of void in the universe. On the one hand, it is impossible for tetrahedra and icosahedra to completely fill a space without leaving gaps,[1] contradicting the absence of void *within* the universe. On the other hand, the shape of the universe, which is composed of polyhedral particles, cannot be a perfect mathematical sphere. This implies the impossibility of the universe's rotation due to the absence of void *outside* it. While the first paradox and the search for a solution have generated considerable scholarly attention, the second paradox seems to have been overlooked by both post-Aristotelian commentators in Antiquity and modern scholars.

The article is divided into four parts. First, we demonstrate that Timaeus' spherical universe cannot be an exact mathematical sphere because it is composed of polyhedra. Then, we provide a detailed demonstration of the mathematical impossibility of Timaeus' universe rotating, as it differs from a perfect mathematical sphere. Moreover, we show that the text itself hints at such a conclusion and that it is highly unlikely Plato was unaware of this issue, challenging the common and simplistic explanations by ignorance or lack of interest. Finally, we propose a solution to this seemingly intractable paradox, based on the fundamental—though highly controversial for modern readers—idea that the universe is a living being with a body and a soul. Ultimately, one must shift from a mathematical perspective to a biological one, taking into account the universe's dynamic nature as '*a single living being that contains within itself all living things, mortal or immortal*'. [2]

The questions addressed here carry significant philosophical implications. These include the relationship between phenomena and principles, particularly the necessity of an interpretation of the *Timaeus* from both philosophical and astronomical perspectives, as noted by Broadie 2012. Additionally, our analysis supports the interpretation that Plato does not reject the sensible universe in favor of mathematics, contrary to the assertions of some scholars.[3] Rather, Plato posits that the entire world follows certain mathematical principles. Whether or not he explicitly stated—as Simplicius reports that Sosigenes claimed—that mathematicians need 'to save the phenomena',[4] this is what is achieved in the *Timaeus* with respect to the primary astronomical knowledge of his time.

These problems are also related to broader philosophical questions. They touch on the relationship between the whole and the parts, which Harte 2002 explores extensively, dedicating an entire chapter to the *Timaeus* (specifically section 4.4, 212-266). They are also connected to the challenging question of the ontology of mathematics, particularly the relationship between mathematics and the intelligible and sensible realms,[5] as illustrated by the intelligible model of the universe.[6] Naturally, these aspects cannot be fully addressed within an article of reasonable length. Nevertheless, they underscore a fundamental point about the *Timaeus*: the interconnection of Plato's philosophical principles and the cosmological construction of the *Timaeus*.

---

[1] We use the usual translation of Aristotle's term 'τόπος' in *De Caelo* III, 8.306b3-10. It is not clearly defined by him; however, it refers probably to some 'simple' regular solid, such as a ball or a cube. For a more detailed discussion of this issue, see Brisson and Ofman 2025, particularly section II.
[2] 'ζῷον ἓν ζῷα ἔχον τὰ πάντα ἐν ἑαυτῷ θνητὰ ἀθάνατά τε.' (69c2)
[3] See for instance the influential books of Thomson 1965 or Farrington 1949. See also Lloyd (1968), p. 80-81.
[4] 'σῴζειν τὰ φαινόμενα' (*In Aristotelis De Caelo Commentaria*, II, 12, 488.14-24); see also Duhem 1908.
[5] E.g. Fronterotta 2007 studies the connection between intelligible ideas and the mathematical order of the universe.
[6] For the important issue of the relation of the model and its image, see for instance Prior 1983; see also Sedley 2016

## I. About the term 'sphere'

In Antiquity, mathematical terminology was not always entirely fixed. In geometry, terms describing figures often exhibited ambiguity, as the same word could refer to both the figure's boundary and to the interior space it enclosed. For plane figures, this meant that a single term might describe a line (the boundary of the surface) and the surface itself, although the latter was more common. For instance, the term 'circle' ('κύκλος') typically referred to what modern mathematics defines as a 'disk', while the term 'circle' is reserved for the curve that forms its circumference.[7]

This is made clear in Euclid's *Elements*, which begins with the following definition of a circle:

> A circle is a plane figure contained by one line such that all the straight lines falling upon it from one point among those lying within the figure are equal to one another (def. I.15).

Similarly, in three-dimensional space, the same term is often used to refer to both a three-dimensional solid and its boundary, a surface, although the primary usage applies to the former. For instance, the sphere ('σφαῖρα') is defined in Euclid's treatise as the solid of revolution generated by rotating a semicircle around its diameter

> When, the diameter of a semicircle remaining fixed, the semicircle is carried round and restored again to the same position from which it began to be moved, the figure so comprehended is a sphere (def. XI.14).

Since the 'sphere' is the solid formed by the revolution of a 'semicircle'—which is a surface, according to Euclid's earlier definition of the 'circle'—around its diameter, the 'sphere' is a three-dimensional solid. However, its boundary is also referred to in the literature as a sphere. In modern mathematics, to avoid such ambiguity, a sphere of radius $R$ and center $O$ consists of all points at a distance $R$ from $O$, while the associated solid, consisting of all points at a distance less than or equal to $R$ from $O$, is called a 'ball'.[8] Thus, the 'sphere' is indeed the boundary of the 'ball'.

Usually, the ambiguity is resolved by the context, which typically concerns one of these two figures. However, here we need to consider both the surface and the solid together. This may cause some difficulty, as both figures are involved in the same passage. To avoid these potential issues, we will use modern terminology unless explicitly stated otherwise. Hence, the 'sphere' will refer to the two-dimensional surface, while the 'ball' will refer to the solid whose boundary is the 'sphere'.

## II. The spherical universe

Some scholars have suggested that the model of the universe might be more like a dodecahedron—the only regular polyhedron not mentioned by Timaeus as a basic particle—than a sphere. They base this claim on a single sentence that alludes to a connection between the universe and the dodecahedron (55d).[9] The dodecahedron is said to be the regular *polyhedron* that most closely approximates the shape of the universe—essentially, the best approximation of the universe's shape among the five regular polyhedra (55c; cf. also *Phaedo* 110b). Hence, it is used by the *demiourgos* to 'paint the star constellations' on the boundary of the universe (55d), probably because it is easier to paint on a flat rather than a curved surface. The claims about the importance of the dodecahedron in shaping the universe seem to originate from interpretations related to its alleged use by Pythagoreans and their influence on Plato.[10] In any case, if the model of the universe were a dodecahedrlon or any polyhedron, the mathematical paradox that we identify in *Timaeus*' cosmology would become even more evident. Thus, we will not consider this point further.

---

[7] A similar ambiguity is found in the *Theaetetus* regarding the term 'powers' ('*dynameis*'). It may refer to the square as a two-dimensional surface or to its perimeter, that is, a one-dimensional line.

[8] Timaeus asserts that the sphere is the solid whose center is equidistant from all its extremities in every direction (33b4-5). Since he refers to its extremity ('μέσης'), it is clear that it has the ball in mind. In any case, this definition is different from Euclid's.

In fact, Timaeus repeatedly asserts that the shape of the universe is spherical (33b4, 44d4, 62d1, 63a5), i.e., that the entire universe is a three-dimensional ball, and its boundary is a two-dimensional sphere. The universe is composed of four types of basic particles—fire, air, water, and earth—each corresponding to one of the four regular polyhedra: the tetrahedron, octahedron, icosahedron, and cube, respectively (56c-57c). It is entirely filled with these particles, which are too small to be visible (56c1-2), with 'no empty *space* left over' ('κενὴν χώραν οὐδεμίαν ἐᾷ λείπεσθαι', 58a7).

This construction aligns with the common ancient Greek conception of the universe, based on these four elements. However, Timaeus introduces something completely new by means of a geometrical feature: that all these basic particles correspond to regular polyhedra. Furthermore, these tiny particles are defined by their faces, which are either equilateral triangles or squares. Both triangles and squares are, in turn, composed of *two* fundamental right triangles: one isosceles, the other half-equilateral. Each face of the cubic particles is formed by four isosceles triangles, whereas each of the other three particles is composed of six half-equilateral triangles, as shown below:

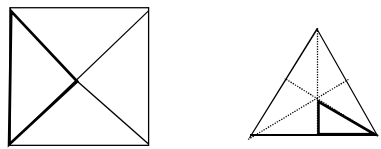

Figure 1

These polyhedral particles are in continuous motion, constantly being destroyed into and recomposed from basic right triangles (isosceles or half-equilateral) (56c-57c). However, using polyhedra to form the ball of the universe presents a difficulty that is studied in Brisson and

---

[9] 'One other construction, a fifth, still remained, and this one the god used for the whole universe, embroidering figures on it.' Some modern scholars, relying on *Timaeus* 55c, consider that the dodecahedron—the only regular polyhedron not used as a basic particle—must be regarded as an essential figure in *Timaeus*' universe and, more generally, in Plato's texts (Petrucci 2021, 2022). This perspective may even lead to the claim that the dodecahedron represents the shape of the universe. A recent paper by theoretical physicists may support this thesis (Luminet 2003). However, this is an absurd assertion within the framework of *Timaeus*, where the dodecahedron's sole function is to serve as the structure upon which the stars are painted on the heavens. Also, see Gregory's note on 55c, though he correctly considers this inconsistent with Timaeus' description of the universe as spherical (Waterfield 2008, 144). As Burnet, quoted by Cornford, explains: 'If the material were not flexible, we should have a regular dodecahedron; as it is flexible, we get a ball.' (Cornford (1937), p. 219). In any case, a dodecahedral universe would render self-rotation impossible, as explained by Aristotle (*De Caelo* II, 4.28a11-23).
[10] See for instance Losee 2001, 16-7, 42-4.

Ofman 2025. In the present article, we consider a different paradox—this time concerning the spherical boundary of the universe, i.e., the two-dimensional sphere that encloses it.

Since the universe is entirely filled with the polyhedral basic particles, its boundary is composed of their faces—hence, of triangles and squares or parts of thereof. Because a mathematical sphere cannot be formed from plane surfaces, just as a mathematical circle cannot be formed from straight lines, the boundary of the universe cannot be a perfect mathematical sphere. It cannot be as smooth as the latter; it must contain hollows and bumps. A very rough cross-section of the sphere of the universe is shown in Figure 2 below:

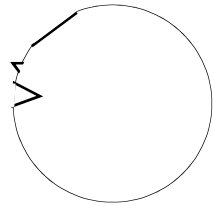

Figure 2

In the above figure:

- the deviations from a perfect mathematical sphere are greatly exaggerated, since they result from the faces of the basic particles, which are themselves so small as to be invisible.
- Only three 'defects' are shown, whereas they are present everywhere along the boundary of the universe.

Because the ratio of the radius of the universe to the size of the basic particles is so immense,[11] the defects, though innumerable, are almost infinitesimal. Its shape is therefore extremely close to being a perfect mathematical sphere—much closer than any sphere produced by a human being—because it is the result of the work of a god. Yet, as close as it is, it is *not* a perfect mathematical sphere.

Timaeus emphasizes the smoothness of the spherical shape of the universe (33b4-c1). However, a little further on, he states that the *demiourgos* shapes the universe like a potter working on a lathe ('ἐτορνεύσατο' 33b5) to give it a spherical form (33b7). This imagery highlights his view of world-shaping as a material task, portraying the *demiourgos* as a craftsman and hinting at the similarity between divine and human works. This suggests that similar difficulties and limitations arise in both cases, differing only in degree. However, the quality of the product is vastly superior when the *demiourgos* is at work compared to a human craftsman.

Specifically, the shape of the spherical universe is *nearly* a mathematical sphere, ensuring that the result is the most beautiful thing 'of all the things that have come to be ('γεγονότων')' (29a5). However, although the *demiourgos* always produces the best and most beautiful outcomes, this is achieved only within the constraints of what is 'as far as possible' (53b5-7; see 29e2-3, 30a1-2, 32b2-7, 32d1-33a1, 37d1-2, 39d8-e2, 89d4-7). This limitation applies even to the astral gods (40b2-4).[12] Consequently, the text aligns with the construction described in the previous section, in which the boundary of the universe is formed from triangles or squares: the shape of the universe is not a perfect mathematical sphere.

---

[11] For discussions of the immense size of the universe as understood in the 4th century BCE, see, for instance, *Stateman* 270a6-8; *De Caelo* II 13.293b32-294a8, 298a15-17;

### III. The paradox

Timaeus' universe is finite, completely filled, and contains everything (33c5-6), with no void either inside or outside it—a view shared by Aristotle. This implies that there is nothing outside the universe, specifically *no space*:

> It is therefore evident that there is also no place or void or time outside the heaven. For in every place body can be present.[13]

There is nothing outside (or beyond) the universe, or more precisely, there is no "outside" (or "beyond") to the universe. This is a difficult issue because it is almost impossible to represent the universe, as all our representations take place within space. In fact, Aristotle, and probably Plato, were perfectly aware of this difficulty and the criticisms their cosmologies would face on this point. Among the discussions of his predecessors, Aristotle finds their objections to a finite universe without an external void particularly challenging.[14] It seems that Parmenides was the first natural philosopher to propose this idea of a finite universe (Furley 1987, 57). Melissus, a philosopher closely aligned with Parmenides and who agreed with most of his philosophy, rejected this particular claim (Sedley 2003, 79). The next known Greek philosopher to support a similar standpoint was Plato, followed by Aristotle.[15] However, contrary to Parmenides, they consider the universe to be in motion. Specifically, Timaeus claims that the whole universe rotates on itself, propelled by the so-called 'circle of the Same', while for some parts this rotation is combined with another, propelled by the so-called 'circle of the Different'. The figure below summarizes his construction (for clarity, only the Earth and the sun are represented):

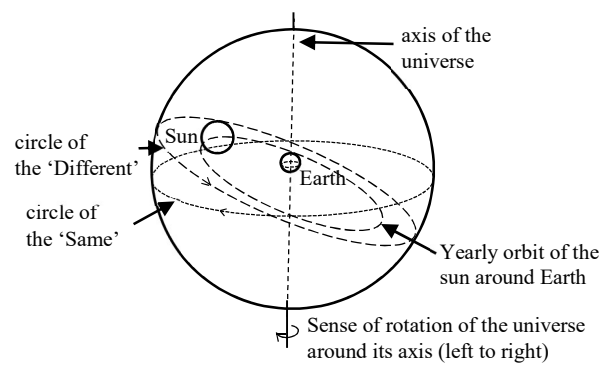

Figure 3

---

[12] Cornford emphasizes that, due to the 'inherent nature' of the materials, the perfection of the *demiourgos*' products is limited (*ibid.*, 37, 41). See also Broadie 2012, 244. On the ontological distinction between mathematical figures and their sensible realizations, see *Republic* VI, 510d5-511a1, 527b5-8.
[13] "Ἅμα δὲ δῆλον ὅτι οὐδὲ τόπος οὐδὲ κενὸν οὐδὲ χρόνος ἐστὶν ἔξω τοῦ οὐρανοῦ. Ἐν ἅπαντι γὰρ τόπῳ δυνατὸν ὑπάρξαι σῶμα·" (*De Caelo* I 9.279a11-12).
[14] Cf. *Physics* III 4.203b15-30. For some discussions on this question, see Sorabji 1988, 125; Alexander, *Quaestiones* 3.12, 106.35-107.4; Archytas of Tarentum quoted in Simplicius, *On Aristotle's Physics* 467.26-32 = DK 47 A 24 (Huffman 2005, 541).
[15] However, Diogenes Laertius (*Lives of Eminent Philosophers*, book IX) reports that Heraclitus considered a finite universe, but since he does not provide his source, it is difficult to ascertain the reliability of this claim.

The universe rotates around its polar axis without changing its position, that is, while all its points, except those on the polar axis change position, the universe as a whole remains at the same place (34a1-4, 37c6).[16] This is possible for a ball due to its symmetry, specifically because a ball is a solid of revolution formed by rotating a half-disk around its diameter. This point is essential because, with no void space outside the universe, it *cannot* change place, as it would, for instance, through translation. This is made clear by Aristotle's claim in the following passage of *De Caelo* that, among all spatial solids, the ball is the only one capable of movement without changing place.[17] Because of its importance, we quote the entire passage:

> Again, since the universe ('τὸ πᾶν') revolves, as observed and in theory, in a circle, and since it has been shown that outside the farthest circumference there is neither void nor place, from these grounds also it will follow necessarily that the universe is spherical. For if it is to be rectilinear in shape, it will follow that there is place and body and void outside the farthest circumference. For a rectilinear figure as it revolves never continues in the same room, but where formerly was body, is now none, and where now is none, body will be in a moment because of the change of position of the corners. Similarly, if the universe had some other figure with unequal radii from the centre, if, for instance, it were lentiform, or oviform, in every case we should have to admit space and void outside the moving body, because the whole body ('τὸ ὅλον') would not always occupy the same space.[18]

Aristotle asserts that the form of the universe cannot be a polyhedron but *must* be a ball. Actually, at the beginning of his account of the universe's formation, Timaeus emphasizes that; among all 'the seven motions', the *demiourgos* granted the universe only one: 'spinning upon itself', or rotation around its axis. The *demiourgos* 'set it turning continuously in the same place', in contrast to the six other types of movement (34a; cf. also *Laws* X, 893b ff). Later, describing the stars as the bodies of gods, Timaeus again emphasizes that one of their two movements is 'rotation, an unvarying movement in the same place' (40a8).

Why is the spherical shape so important in this context? The reason lies in the fact that, unlike rectilinear solids, a sphere can rotate *without changing place*. To better understand this, let us compare the effects of the ball's rotation with those of the cube's rotation, considering for simplicity a cross-section of these figures by a plane:

---

[16] A change of place for a body implies that the body's position in space at two different moments is not the same. In more Aristotelian terms, its occupies two distinct locations in space. Such a change does not necessarily entail a translation. For example, a cube rotating around an axis changes place without undergoing any translation motion (cf. Figure 4, *infra*).

[17] Actually, any solid of revolution, such as a cylinder or a cone, can rotate around its axis of revolution without changing place; however, owing to its strong symmetry, the ball can achieve this around any of its axes, meaning any of its diameters.

[18] ''Ἔτι δὲ ἐπεὶ φαίνεται καὶ ὑπόκειται κύκλῳ περιφέρεσθαι τὸ πᾶν, δέδεικται δ' ὅτι τῆς ἐσχάτης περιφορᾶς οὔτε κενόν ἐστιν ἔξωθεν οὔτε τόπος, ἀνάγκη καὶ διὰ ταῦτα σφαιροειδῆ εἶναι αὐτόν. Εἰ γὰρ ἔσται εὐθύγραμμος, συμβήσεται καὶ τόπον εἶναι ἔξω καὶ σῶμα καὶ κενόν. Κύκλῳ γὰρ στρεφόμενον τὸ εὐθύγραμμον οὐδέποτε τὴν αὐτὴν ἐφέξει χώραν, ἀλλ' ὅπου πρότερον ἦν σῶμα, νῦν οὐκ ἔσται, καὶ οὗ νῦν οὐκ ἔστι, πάλιν ἔσται, διὰ τὴν παράλλαξιν τῶν γωνιῶν. Ὁμοίως δὲ κἂν εἴ τι ἄλλο σχῆμα γένοιτο μὴ ἴσας ἔχον τὰς ἐκ τοῦ μέσου γραμμάς, οἷον φακοειδὲς ἢ ᾠοειδές· ἐν ἅπασι γὰρ συμβήσεται καὶ τόπον ἔξω καὶ κενὸν εἶναι ἧς φορᾶς, διὰ τὸ μὴ τὴν αὐτὴν χώραν κατέχειν τὸ ὅλον. ἧς φορᾶς, διὰ τὸ μὴ τὴν αὐτὴν χώραν κατέχειν τὸ ὅλον.' (II 4.287a11-23; translation from the Ross-Barnes edition, slightlymodified).

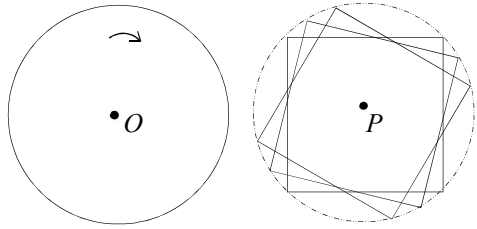

Figure 4

When the disk rotates in the plane, the symmetry of the circle ensures that this does not cause any change in its position; there is no difference of its position whether the circle is rotating or remaining stationary. However, when the square rotates, it changes its position and traces out another figure—the dotted circle in the figure.

This holds true for a perfect *mathematical* ball. However, what happens when dealing with a ball whose shape deviates even slightly from such a perfect mathematical form? Even the slightest irregularity will result in some displacement of the ball in space. Let us again consider a cross-section of a sphere, denoted as **D1** in Figure 5 below, with a single irregularity: a 'peak' at point *A*:

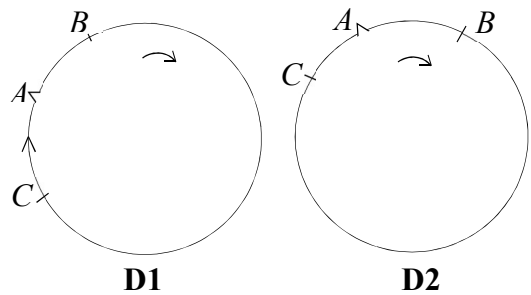

Figure 5

When the ball rotates around its center, **D1** transforms into **D2**, causing the 'peak' at *A* to shift to where *B* was previously, which lies outside of **D1**, as shown in **D2**. In Aristotle's terms, where there was a peak outside the disk at *A* in **D1** there is now nothing outside the disk at *A* in **D2**, and where there was nothing outside the disk at *B* in **D1**, there is now something outside the disk (i.e., the peak) in **D2**.

For both concave irregularity (Figure 6) and plane irregularity (Figure 7), the same situation occurs, as illustrated:

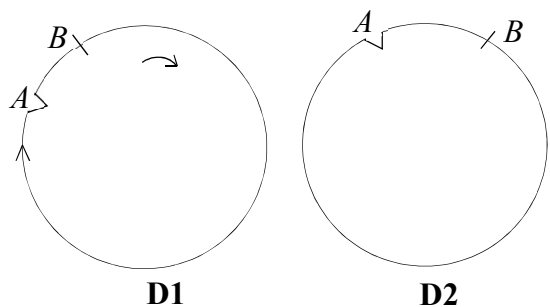

Figure 6

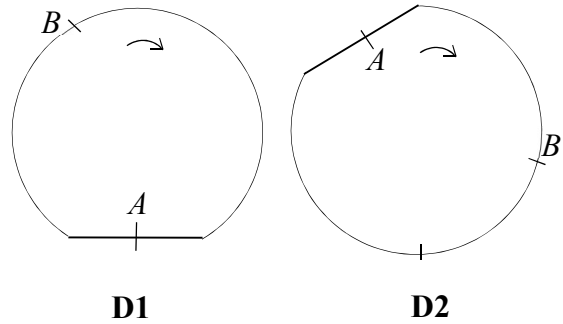

Figure 7

In Figure 6, at point *A* in **D1**, there a 'cavity' that is filled in **D2**. In Figure 7, the 'plane' part at *A* in **D1** becomes circular in **D2**. In Aristotle's terms, in both figures, where there was nothing (the 'cavity' or the 'plane' part) at *A* in **D1**, there is now something (the disk) in **D2**. Conversely, where there was something at *B* (the disk) in **D1** there is now nothing (because of the cavity or the plane part) in **D2**. Consequently, any ball that is not a perfect mathematical ball, when rotating around its axis, causes a change in its position.

Since the boundary of the universe is not a perfect mathematical sphere, any rotation would necessarily imply a change of position. However, any change in position is impossible for the universe because the universe constitutes the totality of the space. More precisely, considering the drawings in the figure as simplified representations of the universe, the 'peak' at *A* would need to move to a 'peak' where *B* previously was, thus outside the limits of **D1**—i.e., outside the universe. Since there is no space outside the universe, nothing can move beyond it—that is, outside **D1**. Thus, it is impossible for the 'peak' at *A* to move (Figure 5). The same applies for a cavity (Figure 6) or a straight part at (Figure 7). In other words, **D1** cannot transition to **D2**.

Consequently, the universe cannot rotate, in contradiction with Timaeus' many claims that the universe is rotating. This is precisely Aristotle's argument in the above-quoted passage of *De Caelo*:

> if the universe had some other figure with unequal radii from the centre, if, for instance, it were lentiform, or oviform, in every case we should have to admit space and void outside the moving body, because the whole body ('τὸ ὅλον') would not always occupy the same space.

Aristotle's statement can be rephrased to suggest that if some points on the boundary of a figure are not equidistant to some central point, it would not be able to rotate. This does not concern only polygonal figures. For example, an ellipse will change its position when rotating:[19]

---

[19] This is likely the reason for Aristotle's error in asserting that only the sphere can move without changing place. He was probably considering the rotation of a figure in the plane, where, indeed, a rotating oviform or lentiform shape does change position.

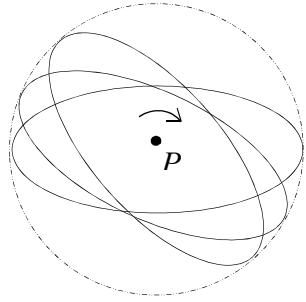

Figure 8

To conclude this section, let us address an objection, that, contrary to our statement, there is no paradox simply because the spherical boundary of the universe can be an exact mathematical sphere. This objection is based on the hypothesis that the basic triangles forming the faces of the particles are not mathematical entities but material bodies, whose shapes do not correspond to perfect mathematical polyhedra. It is thus argued that, since the impossibility applies only to mathematical polyhedra, these imperfect particles can completely fill a ball. Consequently, the boundary formed by the faces of these polyhedra could compose a perfect mathematical sphere, implying that the boundary of the universe could indeed be a perfect mathematical sphere.

We will not revisit the question of whether the basic triangles are mathematical or material, which was addressed in previous articles,[20] nor will we assess the consistency of these claims. We will focus solely on the reasoning itself and argue that it is flawed. First, any imperfections in the faces of these particles would accumulate, compounding with other irregularities previously discussed and thereby increasing the discrepancy between the spherical shape of the universe and a perfect mathematical sphere. Moreover, to form a sphere, each polyhedron would need to have at least one perfectly spherical face. Additionally, all these faces would need to lie on a sphere with the same radius—the radius of the universe itself. Finally, these 'polyhedra' would need to fit together perfectly, which is virtually impossible given that all particles are in continuous motion.

## IV. The usual 'solutions'

As with many difficulties in Plato's texts, some scholars suggest that Plato might have considered the difference between the shape of the 'spherical' universe and a perfect mathematical sphere to be inconsequential. They argue that, since Plato does not address this question, he either dismisses it as unimportant or fails to recognize it as an issue. Others propose that Plato may not have been aware of the problem at all—the impossibility of constructing a sphere from flat surfaces.

First, the claim that the absence of an explicit discussion indicates Plato's disregard for or ignorance of this question is based on a flawed anachronistic assumption: that modern readers of Plato share the same knowledge as his contemporaries.[21] This assumption neglects the fact that many ideas obvious to Athenians in the 4th century BCE may be entirely unfamiliar to contemporary readers. Explanations or clarifications essential for modern audiences may have seemed trivial to Plato's original audience. Therefore, the absence of an explicit explanation does not necessarily imply that Plato considered the issue unimportant or failed to address it.

---

[20] Cf. Brisson and Ofman 2021 and 2022.
[21] As Burnyeat notes in the case of mathematics in Plato's *Republic* (Burnyeat 2000, 24).

Rather, it likely reflects his assumption that his audience would understand it immediately due to the general knowledge of the time.

Yet, examples of such discrepancies are plentiful throughout Plato's works, especially concerning mathematical matters. In fact, most of the main debated issues among Platonists arise from such discrepancies. For instance, there is the problem of the meaning of '*dynameis*' in the famous 'mathematical passage' of *Theaetetus* (147d-e); the interpretation of the second mathematical example in *Meno* (86e-87b); and the question of whether the Earth rotates on its axis in *Timaeus* (40b8-c3), along with Aristotle's commentary on this question (*De Caelo*, II, 13, 293b30-32). All these issues have given rise to extensive discussions and studies in modern times. Furthermore, Platonists such as Proclus and Simplicius, in their commentaries on the *Timaeus* and *De Caelo* respectively, struggled to understand these texts. This demonstrates that, even in late Antiquity, much of the knowledge evident in the 4th century BCE had already been lost.[22] Moreover, in some cases, it is simply impossible to provide all the details concerning complex issues, as completely satisfactory explanations are not always possible, as Plato explains in the passage 68c-d.

The suggestion that Plato disregards the fact that the spherical universe is not a perfect mathematical sphere—because he does not consider small details—aligns with the explanation proposed by some scholars regarding the 'interstices' between the polyhedral basic particles filling the universe. The voids between particles are so minute that Plato may have deliberately disregarded such details.[23] Similarly, they argue that the 'imperfections' relative to a perfect mathematical sphere—i.e., the hollows and bumps on the surface of the universe—are so minute that Plato simply ignores them.

However, this interpretation makes no sense, because the mathematical impossibility does not depend on the size of the defect. Any deviation, no matter how slight, from the shape of a perfect mathematical sphere implies the impossibility of the universe moving. As Aristotle states in *De Caelo* (I 5.271b5 *ff.*), echoing Plato himself (*Cratylus*, 436d), even a small error in mathematics can have immense consequences.[24]

This principle applies to the shape of the universe: although the inaccuracies relative to a perfect mathematical shape are minute, they occur along its entire boundary, rendering them nearly infinite in number. Moreover, nothing, however small, can move outside of the universe—that is, outside of space itself. Since *Timaeus*' universe is made of polyhedra, there are numerous corners along its periphery. When the universe rotates, these corners also rotate, which implies, as Aristotle says, that the universe 'never continues in the same room, but where formerly was body, is now none, and where now is none, body will be in a moment because of the change of position of the corners.' This is precisely Aristotle's main argument against a non-spherical universe, making it highly implausible that Plato would disregard such a fundamental mathematical property.

Plato's supposed ignorance seems even less plausible. A cursory examination of Plato's works—ranging from *Republic* to *Laws*—demonstrates his excellent mathematical skill and adept application of fundamental mathematical principles. Moreover, Plato's profound knowledge across every field he addresses, extends to every field that he discusses, from the biology of the human body in the *Timaeus* to the practicalities of hunting and fishing in the *Sophist* and the *Laws*. Considering the central role mathematics plays in his philosophy, it is highly unlikely that he lacked a comparable level of mathematical understanding.

---

[22] The question of whether the Earth moves or remains immobile is the subject of a forthcoming article.
[23] Cf. O'Brien 1984, 360-362; for a dissenting view, see Brisson and Ofman 2025.
[24] A problem well known by astronomers in creating calendars over the long term. A mere error of one second per day will result in a full day of discrepancy over 3600 years. This is why the lengths of months are irregular (sometimes 28, 29, 30, or 31 days) and why an extra day is added every four years (a leap year).

In the specific case of our problem, is it conceivable that Plato either did not care about or was unaware of the geometrical impossibility of the construction a sphere from linear surfaces? To address this question, one must consider the state of knowledge on such issues in the 4th century BCE. As usual, Aristotle's testimonies are crucial. His primary discussions of these questions are found in several works, most notably *Sophistical Refutations* (11, 171b12-15; 172a2-3), *Prior Analytics* (II, 25.69a25-34), and *Physics* (I, 2.185a15-19), where he examines the famous problem of the so-called 'squaring the circle'.

The problem of 'squaring the circle' involves finding a square with the same area as a given disk. This challenge captured the attention of thinkers in the mid-5th century BCE, including renowned mathematicians such as Hippocrates of Chios and sophists like Bryson and Antiphon (Aristotle, *Physics* I, 2.185a15-19). By the late 5th century, the problem was widely regarded as unsolvable, a view even reflected in popular culture. In Aristophanes' *Birds*, the character Meton is dismissed as a charlatan for claiming he can solve it (v. 990-1020). According to Plutarch, Anaxagoras had already written about this problem.[25]

By Plato's time, it had become an iconic example of an impossible task. The connection between this problem and the 'linearization' of the circle lies in the method of inscribing a regular polygon within the disk, with so many sides that it effectively coincides with the circle—an approach likely used by Antiphon and Bryson.

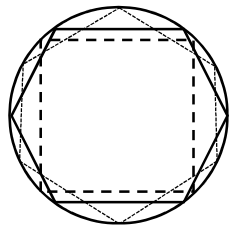

Figure 9

The method begins by inscribing a square within the disk, followed by a regular hexagon, then a regular octagon, and so on. As the number of sides of the inscribed regular polygons increases, their areas grow and progressively approach the area of the disk. Aristotle observes that Antiphon and Bryson committed fallacies in claiming to have squared the disk (*Sophistical Refutations*, 171b16). However, Archimedes later employed a similar method to approximate the ratio of the perimeter to the diameter of a circle—in modern terms, the value of π.[26]

By the 5th century BCE, it was widely accepted that no polygon could completely fill a disk, implying that the circle cannot be constructed from straight lines (cf. Figure 9 above). In three-dimensional space, 'squaring of the circle' extends to 'cubing the ball', which refers to filling the ball entirely with polyhedra, without any gaps or overlaps. This implies that its boundary, the sphere, would be composed of flat figures— such as triangles, squares, or parts of thereof.

---

[25] *On Exile* 17, 607e, f = DK59A8, see Sider 2005.
[26] In his treatise *On the measurement of a circle*, cf. Heath 1949, 48-50.

It is inconceivable that Plato was unaware of or indifferent to such a central problem in the mathematics of his time, especially in light of Aristotle's criticism, that *Timaeus*' basic particles cannot completely fill the universe.[27]

## V. The dynamical solution

The paradox is indeed intractable for a mathematical solid. However, *Timaeus*' universe is a 'living body' ('ζῶν', 37c6, 37d3, …), and must be approached dynamically. As Timaeus emphasizes, all bodies within the universe move along with it. The mechanism of these movements, in a universe without void, is explained later in the passage 79b-80c, which begins with the exposition of the mechanism of respiration in living beings:

> Since there is no void into which anything that is moving could enter, and since the air we breathe out does move out, away from us, it clearly follows that this air doesn't move into a void, but pushes the air next to it out of its place. As this air is pushed out, it drives out the air next to it, and so on, and so inevitably the air, displaced all around, enters the place from which the original air was breathed out and refills that place, following hard on the breath. This all takes place at once ('ἅμα'), like the rotation of a wheel, because there is no such thing as a void.[28]

Some particles of air exiting the body displace others outside the body, which, in turn, displace additional particles, and so forth. This process ensures that new particles of air replace those that have left the body.[29] Timaeus' model likens this process to a wheel. The various parts of a wheel move together simultaneously. When the wheel rotates, the part at a specific position displaces another, setting off a chain reaction, with each part sequentially replacing the previous one. This ensures that no void appears between the parts of the wheel throughout the entire movement. Timaeus further explains that this principle applies not only to respiration but to general movements, including water currents, lightning, and the attractive forces of amber and lodestone. As Timaeus asserts, 'any careful investigator will discover', since there is no void, these things push themselves around into each other; all things move by exchanging places, each to its own place, whether in the process of association or of dissociation.[30]

---

[27] *De Caelo* III 8.306b4-15. More specifically, in this passage, Aristotle asserts that, aside from cubes and octahedra, the other two regular polyhedra used by Timaeus as basic particles—the tetrahedron and the icosahedron—cannot 'fill a space' without leaving voids. In Brisson and Ofman 2025, we address this issue in more detail. Furthermore, it is likely that many of Aristotle's criticisms of Plato's theses were discussed extensively within the Academy during Plato's time; see in particular Vegetti 2000, 450-1.

[28] 'ἐπειδὴ *κενὸν* οὐδέν ἐστιν εἰς ὃ τῶν φερομένων δύναιτ' ἂν εἰσελθεῖν τι, τὸ δὲ πνεῦμα φέρεται παρ' ἡμῶν ἔξω, τὸ μετὰ τοῦτο ἤδη παντὶ δῆλον ὡς οὐκ εἰς *κενόν*, ἀλλὰ τὸ πλησίον ἐκ τῆς ἕδρας ὠθεῖ· τὸ δ' ὠθούμενον ἐξελαύνει τὸ πλησίον ἀεί, καὶ κατὰ ταύτην τὴν ἀνάγκην πᾶν περιελαυνόμενον εἰς τὴν ἕδραν ὅθεν ἐξῆλθεν τὸ πνεῦμα, εἰσιὸν ἐκεῖσε καὶ ἀναπληροῦν αὐτὴν συνέπεται τῷ πνεύματι, καὶ τοῦτο *ἅμα* πᾶν οἷον τροχοῦ περιαγομένου γίγνεται διὰ τὸ *κενὸν* μηδὲν εἶναι.' (79b1-c1).

[29] To avoid any confusion, this does not mean that the universe is breathing in any way; this notion is explicitly denied by Timaeus (32c3-4). The concept of breathing is used here an illustration of how Timaeus asserts the *simultaneity* of two different phenomena: the destruction of some polyhedra and their replacement by others. *Simultaneity* is the key feature required to avoid a void, which would necessarily appear in the universe if these phenomena did not occur *simultaneously*.

These transformations—association ('συγκρινόντων') and dissociation ('διακρινόντων', 58c7)—are the cause of motion (58c), and of most sensations (65c-e), particularly sight (67d-e). Even respiration follows these processes, as it is not a simple exchange of particles between the living being and the surrounding air, but rather requires a cycle of cooling and warming (79d-80a).

The paradigm of respiration offers a solution to the difficulty regarding the rotation of the spherical universe, though it is not a perfect mathematical ball. During the universe's rotation, the polyhedra forming any given 'peak' or other irregularities are either in motion or destroyed within the universe, creating free space in both cases. However, due to the absence of void and the internal pressure of the universe, this space is filled 'at once' by polyhedra newly formed from basic right triangles. This process occurs *simultaneously*, ensuring that the 'peak' remains in the same position, even though it is now composed of new polyhedra.

Indeed, the common model of the universe as a hard, solid ball is unsuitable. Instead of imagining the universe as a rigid, bowling-ball-like structure, it is better to compare it to a balloon, where a deformation in one location causes a corresponding deformation in another place. Of course, because we always represent a body as existing within space, any model that situates the universe within space is inherently inadequate — since the universe constitutes the entirety of space. For contemporary readers, a useful analogy might nevertheless be a wheel with a rubber tire.[31] Its form is mostly cylindrical, with a flat portion—the part in contact with the road—represented as line *AB* in Figure 10:

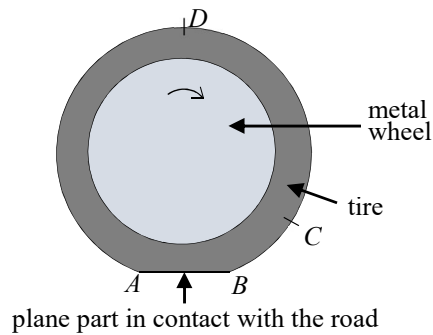

Figure 10

A rigid and solid wheel made of wood or metal with the same plane part in contact with the road, as shown in Figure 11, would be difficult to rotate:

---

[30] 'τὸ δὲ κενὸν εἶναι μηδὲν περιωθεῖν τε αὐτὰ ταῦτα εἰς ἄλληλα, τό τε διακρινόμενα καὶ συγκρινόμενα πρὸς τὴν αὐτῶν διαμειβόμενα ἕδραν ἕκαστα ἰέναι πάντα' (80c3-5). The assertion of the simultaneity of two events, even when one is the cause of the other, may surprise a modern reader accustomed to assuming the necessity of some interval of time between two such events (see, for instance, Gregory 2001, 218, on the decomposition and recomposition of atoms). However, this viewpoint was not necessarily held by the ancient Greeks, as exemplified by the considerations on the freezing of water (Aristotle, *Sense and Sensibilia*, 417b2-3) and the instantaneous perception of light (Aristotle's *On the Soul*, 418b23). For more details, we refer to Brisson and Ofman 2025.

[31] While the universe is spherical and the tire is cylindrical, what is observed is a cross-section—a circular form in both cases, as seen in the drawings. There are other differences: unlike the tire, the universe cannot be inflated or deflated; it is not in contact with anything, as there is nothing outside it, and so on. However, the point here is to illustrate how a non-circular form can rotate without changing position.

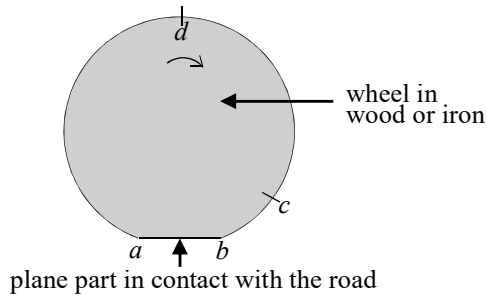

Figure 11

This becomes immediately clear when one tries to drive with a deflated tire. However, with a properly inflated tire, the difficulty disappears because the tire material is supple. From the perspective of overall rotation, the tire's position does not change. This is evident when observed from the viewpoint of someone outside the car but moving alongside it at the same speed, or when the car is placed on a treadmill moving at the same speed as the car, but in the opposite direction, to counteract the translation motion of the wheel. This is depicted in the left drawing of Figure 12.

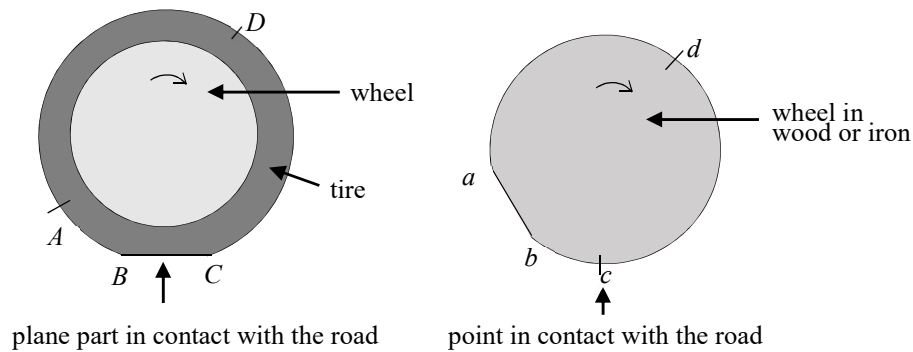

Figure 12

In both the left and right diagrams of Figure 12, all points change positions. However, during its rotation, the overall position of the entire wheel with a rubber tire does not change, as shown in the left diagram. This occurs because the rubber tire, which was compressed by the floor at *A*, is no longer compressed when *A* rotates. This is true for all the points between *A* and *B*, so that *AB* regains its circular form. Meanwhile, the rubber tire, which was not compressed at *C*, becomes compressed on the road. This is true for all the points between *B* and *C*; thus, the flat part at *AB* is replaced by the flat part at *BC*.

This behavior is made possible by the deformability of the supple tire. When a tire rotates, the portion in contact with the road shifts from *AB* to *BC*. The deformed section *AB* moves away and regains its original circular shape, while the part *BC* becomes flat as it comes into contact with the road. This cycle of deformation and recovery repeats continuously as the tire rolls, allowing it to maintain smooth contact with the road while supporting the vehicle's load. In contrast, a solid wheel made of wood or metal changes its overall position during rotation due to its rigidity, as shown in the right-hand diagram of Figure 12.

In other words, the scenario depicted in the right-hand drawing of Figure 7 never occurs for the tire: the 'irregularity' (i.e., the flat part of the tire) remains in the same position as the tire rotates. Hence, even if there were no space around the tire, it would be able to rotate. The same principle applies to the universe: the irregularity stays at the same place while the universe rotates. This principle holds true for any irregularity—whether convex (protruding outward from the spherical form), concave (receding inward), flat, or any combination of these three shapes. To avoid any confusion, just as a rotating wheel is compressed by the road, the universe is compressed by itself due to the absence of exterior space.

Similarly, although every part of the universe is moving and transforming, the universe as a whole neither changes position nor alters its overall shape. While, for the tire, this invariance results from physical laws and holds only to some extent, for the universe, this invariance is an absolute necessity, as it represents the entirety of space. There is no change of position or shape, however minute, of the universe as a whole. This is why it can rotate around its axis, even though it does not have the shape of an exact mathematical sphere.

**Conclusion**.

We have attempted to demonstrate that, when viewed through the lens of conventional cosmological models, *Timaeus*' construction reveals an insurmountable defect rooted in geometry. As Aristotle observes, this difficulty renders *Timaeus*' cosmology contradictory to the most exact of sciences: mathematics. This contradiction arises from two fundamental features of Timaeus' (but also Aristotle's) universe: the absence of void and the rotation of the universe. To grasp the consistency of *Timaeus*' universe, one must shift from a mathematical perspective to a biological one. This implies considering it as a dynamical object, that is a 'living being' rather than a rigid geometrical solid, since it is, in its entirety, in constant motion, with the particles composing it undergoing continuous destruction and reconstruction. We attempted to demonstrate that this method both allows us to highlight a paradox that was almost completely overlooked by Plato's commentators and provides a solution to it.

This dynamical aspect of *Timaeus*' universe has frequently been criticized as inconsistent by many commentators, both ancient and modern, who have analyzed the text through the lens of Aristotle's cosmology. Our primary objective has been to understand *Timaeus*' cosmology using Platonic concepts, demonstrating its consistency from both a philosophical and a physical standpoint, while addressing the criticisms of Aristotle and most modern commentators.